
\def\LaTeX{%
  \let\Begin\begin
  \let\End\end
  \let\salta\relax
  \let\finqui\relax
  \let\futuro\relax}

\def\UK{\def\our{our}\let\sz s}
\def\USA{\def\our{or}\let\sz z}

\UK



\LaTeX

\USA


\documentclass[twoside,12pt]{article}
\setlength{\textheight}{23.5cm}
\setlength{\textwidth}{16cm}
\setlength{\oddsidemargin}{2mm}
\setlength{\evensidemargin}{2mm}
\setlength{\topmargin}{-15mm}
\parskip2mm

\salta


\usepackage[usenames,dvipsnames]{color}
\usepackage{amsmath}
\usepackage{amsthm}
\usepackage{amssymb}
\usepackage[mathcal]{euscript}
%
%
%
%
%
%
%


\definecolor{viola}{rgb}{0.3,0,0.7}
\definecolor{ciclamino}{rgb}{0.5,0,0.5}

\def\betti #1{{\color{Green}#1}}
\def\gabri #1{{\color{viola}#1}}
\def\gianni #1{{\color{blue}#1}}
\def\pier #1{{\color{red}#1}}
\def\gianni #1{#1}
\def\pier #1{#1}
\def\betti #1{#1}
\def\gabri#1{#1}


\bibliographystyle{plain}


%

\finqui

\def\Beq{\Begin{equation}}
\def\Eeq{\End{equation}}
\def\Bsist{\Begin{eqnarray}}
\def\Esist{\End{eqnarray}}

\def\Bthm{\Begin{theorem}}
\def\Ethm{\End{theorem}}

\def\Bprop{\Begin{proposition}}
\def\Eprop{\End{proposition}}
\def\Bcor{\Begin{corollary}}
\def\Ecor{\End{corollary}}

\def\Bdim{\Begin{proof}}
\def\Edim{\End{proof}}
\def\Bcenter{\Begin{center}}
\def\Ecenter{\End{center}}
\let\non\nonumber




\def\step #1 \par{\medskip\noindent{\bf #1.}\quad}


\def\Lip{Lip\-schitz}
\def\Holder{H\"older}
\def\Frechet{Fr\'echet}
\def\aand{\quad\hbox{and}\quad}

\def\lhs{left-hand side}
\def\rhs{right-hand side}


\def\lineariz{lineari\sz}

\def\regulariz{regulari\sz}

\def\bhv{behavi\our}


\def\multibold #1{\def\arg{#1}%
  \ifx\arg\pto \let\next\relax
  \else
  \def\next{\expandafter
    \def\csname #1#1#1\endcsname{{\bf #1}}%
    \multibold}%
  \fi \next}

\def\pto{.}

\def\multical #1{\def\arg{#1}%
  \ifx\arg\pto \let\next\relax
  \else
  \def\next{\expandafter
    \def\csname cal#1\endcsname{{\cal #1}}%
    \multical}%
  \fi \next}


\def\multimathop #1 {\def\arg{#1}%
  \ifx\arg\pto \let\next\relax
  \else
  \def\next{\expandafter
    \def\csname #1\endcsname{\mathop{\rm #1}\nolimits}%
    \multimathop}%
  \fi \next}

\multibold
qwertyuiopasdfghjklzxcvbnmQWERTYUIOPASDFGHJKLZXCVBNM.

\multical
QWERTYUIOPASDFGHJKLZXCVBNM.

\multimathop
dist div dom meas sign supp .


\def\accorpa #1#2{\eqref{#1}--\eqref{#2}}
\def\Accorpa #1#2 #3 {\gdef #1{\eqref{#2}--\eqref{#3}}%
  \wlog{}\wlog{\string #1 -> #2 - #3}\wlog{}}


\def\graffe #1{\mathopen\{#1\mathclose\}}

\def\<#1>{\mathopen\langle #1\mathclose\rangle}
\def\norma #1{\mathopen \| #1\mathclose \|}

\def\iot {\int_0^t}

\def\intQt{\int_{Q_t}}
\def\intQ{\int_Q}
\def\iO{\int_\Omega}

\def\dt{\partial_t}
\def\dn{\partial_n}

\def\cpto{\,\cdot\,}

\def\checkmmode #1{\relax\ifmmode\hbox{#1}\else{#1}\fi}
\def\aeO{\checkmmode{a.e.\ in~$\Omega$}}
\def\aeQ{\checkmmode{a.e.\ in~$Q$}}

\def\aeS{\checkmmode{a.e.\ on~$\Sigma$}}

\def\aaQ{\checkmmode{for a.a.~$(t,x)\in Q$}}


\def\erre{{\mathbb{R}}}




\def\genspazio #1#2#3#4#5{#1^{#2}(#5,#4;#3)}
\def\spazio #1#2#3{\genspazio {#1}{#2}{#3}T0}

\def\L {\spazio L}
\def\H {\spazio H}
\def\W {\spazio W}

\def\C #1#2{C^{#1}([0,T];#2)}


\def\Lx #1{L^{#1}(\Omega)}
\def\Hx #1{H^{#1}(\Omega)}

\def\Cx #1{C^{#1}(\overline\Omega)}

\def\LQ #1{L^{#1}(Q)}

\def\CQ #1{C^{#1}(\overline Q)}

\def\Ldue{\Lx 2}
\def\Linfty{\Lx\infty}

\def\Huno{\Hx 1}
\def\Hdue{\Hx 2}
\def\Hunoz{{H^1_0(\Omega)}}


\def\LQ #1{L^{#1}(Q)}


\let\theta\vartheta
\let\eps\varepsilon
\let\phi\varphi

\let\TeXchi\chi                         
\newbox\chibox
\setbox0 \hbox{\mathsurround0pt $\TeXchi$}
\setbox\chibox \hbox{\raise\dp0 \box 0 }
\def\chi{\copy\chibox}



\def\thetaz{\theta_0}
\def\phiz{\phi_0}
\def\thetaQ{\theta_{\!Q}}

\def\Vz{V_0}
\def\umin{u_{\rm min}}
\def\umax{u_{\rm max}}
\def\phimin{\phi_\bullet}
\def\phimax{\phi^\bullet}

\def\Uad{\calU_{ad}}
\def\uopt{u^*}
\def\thetaopt{\theta^*}
\def\phiopt{\phi^*}
\def\ubar{\overline u}
\def\thetabar{\overline\theta}
\def\phibar{\overline\phi}
\def\thetah{\theta^h}
\def\phih{\phi^h}
\def\zetah{\zeta^h}
\def\etah{\eta^h}
\def\ptilde{\tilde p}
\def\qtilde{\tilde q}

\def\un{u_n}
\def\phin{\phi_n}
\def\thetan{\theta_n}
\def\xin{\xi_n}

\def\normaV #1{\norma{#1}_V}

\let\hat\widehat
\def\Beta{\hat{\vphantom t\smash\beta\mskip2mu}\mskip-1mu}
\def\betaeps{\beta_\eps}

\def\Betap{\Beta_{\eps,p}}
\def\phieps{\phi_\eps}
\def\xieps{\xi_\eps}

\def\betaz{\beta^\circ}

\def\Pi{\hat\pi}

\def\xieps{\xi_\eps}


\def\cz{C_0}

\Begin{document}

\thispagestyle{empty}

\title{\pier{Optimal control} for a phase field system\\[0.3cm]
  with a possibly singular potential%
}

\author{}
\date{}
\maketitle
\Bcenter
\vskip-2cm
{\large\sc Pierluigi Colli$^{(1)}$}\\
{\normalsize e-mail: {\tt pierluigi.colli@unipv.it}}\\[.25cm]
{\large\sc Gianni Gilardi$^{(1)}$}\\
{\normalsize e-mail: {\tt \gabri{gianni.gilardi}@unipv.it}}\\[.25cm]
{\large\sc Gabriela Marinoschi$^{(2)}$}\\
{\normalsize e-mail: {\tt gabriela.marinoschi@acad.ro}}\\[.25cm]
{\large\sc Elisabetta Rocca$^{(3)}$}\\
{\normalsize e-mail: {\tt rocca@wias-berlin.de} \betti{and {\tt elisabetta.rocca@unimi.it}}}\\[.45cm]
$^{(1)}$
{\small Dipartimento di Matematica ``F. Casorati'', Universit\`a di Pavia}\\
{\small via Ferrata 1, 27100 Pavia, Italy}\\[.2cm]
$^{(2)}$
{\small ``Gheorghe Mihoc-Caius Iacob'' Institute of Mathematical Statistics\\ 
\gabri{{}\ \ and Applied Mathematics of the Romanian Academy (ISMMA)}\\
Calea 13 Septembrie 13, 050711 Bucharest, Romania}\\[.2cm]
$^{(3)}$
{\small Weierstrass Institute for Applied Analysis and Stochastics}\\
{\small Mohrenstrasse 39, 10117 Berlin, Germany\\
 \betti{and} \\
\betti{Dipartimento di Matematica, Universit\`a di Milano}\\
\betti{via Saldini 50, 20133 Milano, Italy}}\\
[1cm]
\Ecenter

\Begin{abstract}
\betti{In this paper we study a \pier{distributed} control problem for a phase field system of Caginalp type with 
logarithmic potential. The main aim of this work would be to force the location of the diffuse interface  
to be as close as possible to a prescribed set. However, due to the discontinuous \pier{character} of the
cost functional, we have to approximate it by a regular one and, in this case, we solve the associated control problem and derive the \pier{related} first order necessary optimality conditions.}
\vskip3mm

\noindent {\bf Key words:}
Phase field system, phase transition,
singular potentials, optimal control, optimality conditions,
adjoint state system.
\vskip3mm
\noindent {\bf AMS (MOS) Subject Classification:} \betti{80A22, \pier{35K55,} 49J20, 49K20}.
\End{abstract}

\salta

\pagestyle{myheadings}
\newcommand\testopari{\sc Colli \ --- \ Gilardi \ --- \ Marinoschi \ --- \ Rocca}
\newcommand\testodispari{\sc \pier{Optimal control} for a phase field system}
\markboth{\testodispari}{\testopari}

\finqui
%
%
\vfill

\section{Introduction}
\label{Intro}
\setcounter{equation}{0}

\betti{This paper is concerned with the study of a distributed control problem for a 
Caginalp type PDE system (cf.~\cite{Cag} and \cite{BrokSpr})}
\Beq
  \dt\theta - \Delta\theta + \ell \dt\phi = \sigma
  \aand
  \dt\phi - \Delta\phi + \calW'(\phi) = \betti{\ell}\theta
  \quad \hbox{in $Q := (0,T)\times\Omega $}
  \label{caginalp}
\Eeq
where $\Omega$ is the domain where the evolution takes place,
$T$~is some final time,
$\theta$~denotes the relative temperature around some critical value
that is taken to be $0$ without loss of generality,
and $\phi$ is the order parameter.
Moreover, $\ell$~is a positive coefficient that is proportional
to the latent heat,
and $\sigma$ is some source term. 
Finally, $\calW'$ represent\gabri{s} the derivative of a double-well potential~$\calW$,
and the typical example is
the classical regular potential $\calW_{reg}$ defined~by
\Beq
  \calW_{reg}(r) = \frac 14 \, (r^2-1)^2 \,,
  \quad r \in \erre .
  \label{regpot}
\Eeq
However, different choices of $\calW$ are possible,
and a thermodynamically significant example is given by the so-called logarithmic double-well potential,
namely
\Beq
  \calW_{log}(r) = ((1+r)\ln (1+r)+(1-r)\ln (1-r)) - c r^2 \,,
  \quad r \in (-1,1)
  \label{logpot}
\Eeq
where $c>0$ is large enough in order to kill convexity. 
More generally, the potential $\calW$ could be just the sum $\calW=\Beta+\Pi$,
where $\Beta$ is a convex function that is allowed to take the value~$+\infty$ \pier{somewhere},
and $\Pi$ is a smooth perturbation (not necessarily concave).
In such a case, $\Beta$~is supposed to be proper and lower semicontinuous
so that its subdifferential is well-defined and can replace the derivative
which might not exist. \betti{A typical example is the so-called double obstacle potential 
\Beq
\calW_{obs}(r)=I_{[-1,1]}(r)-cr^2\,,
\label{obspot}
\Eeq
where $I_{[-1,1]}$ denotes the indicator function of the set $[-1,1]$ which takes value 0 in $[-1,1]$ and $+\infty$ outside.
}
Of course, the second equation \eqref{caginalp} becomes a differential inclusion.

\betti{The mathematical literature on \eqref{caginalp} is rather vast and we confine ourselves to quote the pioneering paper \cite{EllZheng} and the more recent ones \pier{\cite{Lau}, \cite{GraPetSch}, \cite{KenmNiez}} dealing respectively with the cases of regular, singular and non-smooth potentials.
}

Moreover, initial conditions like $\theta(0)=\thetaz$ and $\phi(0)=\phiz$ 
and suitable boundary conditions must complement the above equations.
As far as the latter are concerned, 
we take the homogeneous Dirichlet and Neumann boundary conditions, respectively, that~is
\Beq
  \theta = 0
  \aand
  \dn\phi = 0 
  \quad \hbox{on $\Sigma := (0,T)\times\Gamma $} 
  \non
\Eeq
where $\Gamma$ is the boundary of~$\Omega$ and $\dn$ is the (say, outward) normal derivative.
We note that the latter is very common in the literature
and that the former could be replaced by an inhomogeneous~one.

The aim of this paper is to study a related optimal control problem,
the control being associated to the forcing term $\sigma$ 
that appears on the \rhs\ of the first equation~\eqref{caginalp}.
Namely, we take
$\sigma(t,x) = m(x) \, u(t,x)$,
where $m$ is a given nonnegative function on $\Omega$
and $u$ is the control.
Thus, the state system takes the following form
\Bsist
  & \dt\theta - \Delta\theta + \ell \dt\phi = mu
  & \quad \hbox{in $Q$}
  \label{Iprima}
  \\
  & \dt\phi - \Delta\phi + \beta(\phi) + \pi(\phi) \ni \betti{\ell}\theta
  & \quad \hbox{in $Q$}
  \label{Iseconda}
  \\
  & \theta = 0 
  \aand
  \dn\phi = 0
  & \quad \hbox{on $\Sigma$}
  \label{Ibc}
  \\
  & \theta(0) = \thetaz
  \aand
  \phi(0) = \phiz
  & \quad \hbox{on $\Omega$} 
  \label{Icauchy}
\Esist
\Accorpa\Ipbl Iprima Icauchy
and the control $u$ is supposed to vary in some control box~$\Uad$.
We would like to force the location of the diffuse interface of~$\phi$, 
\betti{i.e., of the set $E_\eps(\phi)$ 
where the state $\phi$ takes values between $-\eps$ and~$\eps$\pier{,}
for some given~$\eps>0$}\pier{,} to be as close as possible
to a prescribed set $E\subset Q$.
To do that, by denoting by $\chi_E$ the characteristic function of~$E$
and by $g$ the characteristic function of the interval~$[-\eps,\eps]$,
we introduce the cost functional
\Beq
  \calJ_0(u)
  := \frac 12 \intQ (g(\phi) - \chi_E)^2
  \label{IcostE}
\Eeq
where $(\theta,\phi)$ is the state corresponding to~$u$.
More generally, we could take,~e.g., 
\Beq
  \calJ(u)
  := \frac 12 \intQ (g(\phi) - \chi_E)^2
  + \frac\kappa 2 \intQ (\theta - \thetaQ)^2
  \label{Icost}
\Eeq
where the desired temperature $\thetaQ\in\LQ2$ 
and the constant $\kappa\geq0$ are given.
In this case, the optimal control (if~it exists)
balances the \gianni{closeness} of $E_\eps(\phi)$ to $E$
and the smallness of the difference $|\theta-\thetaQ|$,
depending on the value of the coefficient~$\kappa$.
However, such problems look difficult for every reasonable control box~$\Uad$.
As this is mainly due to the discontinuous character of~$g$,
we replace the characteristic function $g$ 
\gabri{by} 
a continuous approximation of~it
(still denoted by~$g$),
and a possible choice is the following
\Beq
  g(r) := \frac \lambda {((r^2-\eps^2)^+)^2 + \lambda}
  \quad \hbox{for $r\in\erre $}
  \non
\Eeq
where $\lambda>0$ is small.
At this point, we can generalize the problem
and allow $g$ to be any continuous function on $\erre$
satisfying some growth condition that makes the cost functional
meaningful for every admissible control~$u$,
and boundedness is surely suitable.
Moreover, even $\chi_E$ can be replaced by a more general given function.

Thus, the control problem we address in this paper consists in minimizing
the cost functional 
\Beq
  \calJ(u)
  := \frac 12 \intQ (g(\phi) - \chi)^2
  + \frac\kappa 2 \intQ (\theta - \thetaQ)^2
  \label{IdefJ}
\Eeq
depending on the state variables $\theta$ and $\phi$
satisfying the above state system, over all the controls belonging to some control box~$\Uad$,
where $\chi$ and $\thetaQ$ are given in~$\LQ2$, $\kappa$~is a nonnegative constant
and $g$ is a prescribed real function on~$\erre$, that we assume to be at least continuous and bounded.
As far as the control box in concerned, we take
\Beq
  \Uad := 
  \bigl\{ u \in \LQ2 : \ \umin\leq u\leq\umax\ \aeQ \bigr\}
  \label{Iuad}
\Eeq
where $\umin$ and $\umax$ are given bounded functions.

\betti{Let us mention here that \pier{in our approach the existence of an optimal control is proven for a quite general class of potentials $\calW$: indeed, $\calW$} is assumed to be a 
 smooth perturbation of a convex function $\widehat\beta$ possibly taking \pier{the} value $+\infty$ \pier{somewhere}. Notice that all the three examples \eqref{regpot}, \eqref{logpot}, and \eqref{obspot} fit these assumptions. However, we point out that the derivation of the first order necessary optimality conditions can be made only in case of regular (e.g.~\eqref{regpot}) and singular (e.g.~\eqref{logpot}) potentials (cf.\ Section~\ref{OPTIMUM}).
Hence, the main novelty of the present contribution consists in the fact that we can deal with quite general potentials $\calW$ (even singular) in the phase equation and quite general cost functions $\calJ$. Up to our knowledge, indeed, the literature on optimal control for Caginalp type phase field models is quite poor and often restricted to the case of regular potentials\pier{, or dealing with approximating problems when first order optimality conditions are discussed.  In this framework, let us quote the papers \cite{HoffJiang, HKKY} and references therein and also \cite{BBCG, BCF, CGPS, CGS, CMR, LK, SY, SprZheng} for different types of phase field models.}} 

The paper is organized as follows.
In the next section, we list our assumptions, state the problem in a precise form
and present our results.
The well-posedness of the state system
and the existence of an optim\gabri{al} control will be shown 
in Sections~\ref{STATE} and~\ref{OPTIMUM}, respectively,
while the rest of the paper is devoted 
to the derivation of first order necessary conditions for optimality.
The final result will be proved in Section~\ref{OPTIMALITY} 
and it is prepared in Sections~\ref{FRECHET}, \pier{which is devoted to the study of}
the control-to-state mapping.


\section{Statement of the problem and results}
\label{STATEMENT}
\setcounter{equation}{0}

In this section, we describe the problem under \pier{investigation}
and present our results. \pier{From now on, for simplicity and 
without any loss of generality we take}  \betti{$\ell=1$ in  \Ipbl .}
As in the Introduction,
$\Omega$~is the body where the evolution takes place.
We assume $\Omega\subset\erre^3$
to~be open, bounded, connected\betti{, of class $C^{1,1}$,} 
and we write $|\Omega|$ for its Lebesgue measure.
Moreover, $\Gamma$ and $\dn$ still stand for
the boundary of~$\Omega$ and the outward normal derivative, respectively.
Given a finite final time~$T>0$,
we set for convenience
\Bsist
  && Q_t := (0,t) \times \Omega 
  \aand
  \Sigma_t := (0,t) \times \Gamma 
  \quad \hbox{for every $t\in(0,T]$}
  \label{defQtSt}
  \\
  && Q := Q_T \,,
  \aand
  \Sigma := \Sigma_T \,.
  \label{defQS}
\Esist
\Accorpa\defQeS defQtSt defQS
Now, we specify the assumptions on the structure of our system.
We assume that
\Bsist
  & 
  m \in \Lx\infty
  \aand
  m \geq 0
  \quad \aeO
  \label{hplm}
  \\
  & \Beta : \erre \to [0,+\infty]
  \quad \hbox{is convex, proper and l.s.c.}
  \quad \hbox{with} \quad
  \Beta(0) = 0
  \label{hpBeta}
  \\
  & \Pi: \erre \to \erre
  \quad \hbox{is a $C^1$ function}
  \aand
  {\Pi\,}'
  \quad \hbox{is \Lip\ continuous}\betti{\,.}
  \label{hpPi}
\Esist
\Accorpa\HPstruttura hplm hpPi
We set for convenience
\Beq
  \beta := \partial\Beta 
  \aand
  \pi := {\Pi\,}'
  \label{defbetapi}
\Eeq
and denote by $D(\beta)$ and $D(\Beta)$ 
the effective domains of $\beta$ and~$\Beta$, respectively.
Moreover, $\betaeps$~is the Yosida regularization of $\beta$ at level~$\eps$
and $\betaz(r)$~denotes the element of $\beta(r)$ having minimum modulus
for every $r\in D(\beta)$
(see, e.g., \cite[p.~28]{Brezis}).
It is well known that
both $\beta$ and $\betaeps$ are maximal monotone operators
and that $\betaeps$ is even single-valued and \Lip\ continuous. 
Furthermore (see, e.g., \cite[Prop.~2.6, p.~28]{Brezis}),
we~have
\Beq
  |\betaeps(r)| \leq |\betaz(r)|
  \aand
  \betaeps(r)\to\betaz(r)
  \quad \hbox{for $r\in D(\beta)$}.
  \label{propYosida}
\Eeq

Next, in order to simplify notations, we~set
\Beq
  V := \Huno, \quad
  \Vz := \Hunoz, \quad
  H := \Ldue, \quad
  W := \graffe{v\in\Hx2: \dn v=0}
  \label{defspazi}
\Eeq
and endow these spaces with their natural norms.
The symbol $\norma\cpto_X$ stands for the norm in the generic Banach space~$X$,
while $\norma\cpto_p$ is the usual norm in both
$\Lx p$ and $\LQ p$, for $1\leq p\leq\infty$\betti{.}
Finally, for \pier{$v\in\L2X$ the function $1*v$ is defined by 
\Beq
  (1*v)(t) := \iot v(s) \, ds
  \quad \hbox{for $t\in[0,T]$}
  \label{defconv}
\Eeq
(note that the symbol $*$ is usually employed for convolution products).}

At this point, we describe the state system.
Given $\thetaz$ and $\phiz$ such that \pier{%
\Bsist
  & \thetaz \in \Vz
    \label{hpz}&
  \\
  & 
  \phiz \in V
  \quad \hbox{and} \quad
  \Beta(\phiz)\in L^1(\Omega) &
  \label{hpzbis}
\Esist
}%
\Accorpa\HPdati hpz hpzbis
we look for a triplet $(\theta,\phi,\xi)$ satisfying
\Bsist
  & \theta \in \H1H \cap \L\infty\Vz \cap \L2\Hdue
  \label{regtheta}
  \\
  & \phi \in {}\pier{\H1H \cap \L\infty V \cap \L2W}
  \label{regphi}
  \\
  & \xi \in \pier{\L2 H}
  \label{regxi}
  \\
  & \dt\theta - \Delta\theta + \dt\phi = mu
  & \quad \aeQ
  \label{prima}
  \\
  & \dt\phi - \Delta\phi + \xi + \pi(\phi) = \theta
  \aand
  \xi \in \beta(\phi)
  & \quad \aeQ
  \label{seconda}
  \\
  & \theta(0) = \thetaz
  \aand
  \phi(0) = \phiz
  & \quad \aeO .
  \label{cauchy}
\Esist
\Accorpa\Regsoluz regtheta regxi
\Accorpa\Pbl prima cauchy
\Accorpa\Tuttopbl regtheta cauchy

\betti{Our first result, whose proof is \pier{sketched} in Section~\ref{STATE},} 
ensures well-posedness with the prescribed regularity, stability
and continuous dependence on the control variable in suitable topologies.

\Bthm
\label{Wellposedness}
Assume \HPstruttura\ and \HPdati.
Then, for every $u\in\LQ2$, problem \Pbl\ has a unique solution $(\theta,\phi,\xi)$
satisfying \Regsoluz, 
and the estimate
\Bsist
  && \norma\theta_{\H1H \cap \L\infty\Vz \cap \L2\Hdue}
   \qquad
  \non
  \\
  && \qquad {}
   + \norma\phi_{\pier{\H1H \cap \L\infty V \cap \L2W}}
  + \norma\xi_{\pier{\L2 H}}
  \ \leq \ C_1
  \label{stimasoluz}
\Esist
holds true for some constant $C_1$ that depends only on $\Omega$, $T$, 
the structure \HPstruttura\ of the system,
the norms of the initial data associated to~\HPdati\ and $\norma u_2$.
Moreover, if $u_i\in\LQ2$, $i=1,2$, are given
and $(\theta_i,\phi_i,\xi_i)$ are the corresponding solutions,
then the estimate
\Bsist
  && \norma{\theta_1-\theta_2}_{\gabri{\L2H}}
  + \norma{(1*\theta_1)-(1*\theta_2)}_{\L\infty\Vz}
  \non
  \\[1mm]
  && \quad {}
  + \norma{\phi_1-\phi_2}_{\L\infty H\cap\L2V}
  \non
  \\[1mm]
  && \leq C' \, \norma{(1*u_1)-(1*u_2)}_{\L2H}
  \leq C'' \, \norma{u_1-u_2}_{\L2H}
  \label{contdep}
\Esist
holds true with constants $C'$ and $C''$ that depend only on
$\Omega$, $T$, $\pi$ and~$m$.
\Ethm

Some further regularity of the solution is stated in the next \pier{result, whose proof} 
\betti{is given in Section~\ref{STATE}}.

\Bthm
\label{Regularity}
The following properties hold true.

$i)$~\pier{Assume \HPstruttura\ and \HPdati. Moreover, let 
\Beq
\label{p1} 
\gianni{%
  \phiz \in W
  \aand
  \betaz(\phiz) \in H \,.
}
  \Eeq
Then, the unique solution $(\theta,\phi,\xi)$ given 
by Theorem~\ref{Wellposedness} also satisfies
\Bsist
  & \phi \in \W{1,\infty}H \cap \H1V \cap \L\infty W &
  \label{p2}
  \\
  & \xi \in \L\infty H &
  \label{p3}
  \\
   &
  \norma\phi_{\W{1,\infty}H \cap \H1V \cap \L\infty W}
  + \norma\xi_{\L\infty H}
  \ \leq \ C_2
  \label{p4}  
  \\
  &\phi \in \CQ0
  \aand
  \norma\phi_\infty \leq C_2 &
  \label{phicont}
\Esist
for some constant $C_2$ that that depends only on $\Omega$, $T$, 
the structure \HPstruttura\ of the system,
the norms of the initial data associated to~\HPdati , \eqref{p1} and $\norma u_2$.}

$ii)$~If in addition $\thetaz\in\Linfty$ \pier{and $u\in\L\infty H$}, we also have
\Beq
  \theta \in \LQ\infty
  \aand
  \norma\theta_\infty \leq C_3
  \label{thetabdd}
\Eeq
with a similar constant $C_3$ that depends on $\norma\thetaz_\infty$ \pier{and $\norma u_{\L\infty H}$} as well.

$iii)$~By \pier{further assuming $\betaz(\phiz)\in\Linfty$},
we have that $\xi\in\LQ\infty$ and
\Beq
  \norma\xi_{\LQ\infty} \leq C_4
  \label{stimaxi}
\Eeq
with a constant $C_4$ that depends on \pier{$C_3$} and $\norma{\betaz(\phiz)}_\infty$ in addition.
\Ethm

The well-posedness result for problem \Pbl\ given by Theorem~\ref{Wellposedness} 
allows us to introduce the control-to-state mapping~$\calS$ 
and to address the corresponding control problem.
We define
\begin{align}
  & \calX := \LQ\infty
  \label{defX}
  \\
  & \calY := \calY_1 \times \calY_2
  \quad \hbox{where}
  \label{defY}
  \\
  & \calY_1 := \graffe{v\in\LQ2: \ 1*v\in\L2\Vz}
  \label{defYuno}
  \\
  & \calY_2 := \L\infty H \cap \L2V
  \label{defYdue}
  \\
  & \calS : \calX \to \calY ,
  \quad 
  u \mapsto \calS(u) =: (\theta,\phi)
  \quad \hbox{where}
  \non
  \\
  & \quad \hbox{$(\theta,\phi,\xi)$ is the unique solution 
  to \Tuttopbl\ corresponding to $u$}.
  \label{defS}
\end{align}
Next, in order to introduce the control box and the cost functional,
we assume that
\Bsist
  & \umin,\umax\in\LQ\infty
  \quad \hbox{satisfy} \quad
  \umin \pier{{}\leq{}}\umax
  \quad \aeQ 
  \label{hpUad}
  \\
  & g: \erre \to \erre
  \quad \hbox{is continuous and bounded}
  \label{hpg}
  \\
  &  \kappa \in [0,+\infty)
  \aand
  \chi , \, \thetaQ \in \LQ2
  \label{hpJ}
\Esist
\Accorpa\StrutturacCP hpUad hpJ
and define $\Uad$ and $\calJ$
according to the Introduction.
Namely, we~set
\Bsist
  && \Uad := 
  \bigl\{ u \in \calX : \ \umin\leq u\leq\umax\ \aeQ \bigr\}
  \phantom \int
  \label{defUad}
  \\[0.2cm]
  && \calJ := \pier{\calF} \circ \calS : \calX \to \erre
  \quad \hbox{where} \quad
  \pier{\calF} : \calY \to \erre
  \quad \hbox{is defined by}
  \non 
  \\
  && \qquad\qquad  \pier{\calF}(\theta,\phi)
  := \frac 12 \intQ (g(\phi) - \chi)^2
  + \frac\kappa 2 \intQ (\theta - \thetaQ)^2 .
  \label{defI}
\Esist
\Accorpa\Defcontrol defX defI


Here is our first result on the control problem\pier{; for the proof we refer to} \betti{Section~\ref{OPTIMUM}}.

\Bthm
\label{Optimum}
Assume \HPstruttura\ and \HPdati,
and let $\Uad$ and $\calJ$ be defined by \pier{\accorpa{defUad}{defI}}.
Then, there exists $\uopt\in\Uad$ such~that
\Beq
  \calJ(\uopt)
  \leq \calJ(u)
  \quad \hbox{for every $u\in\Uad$} .
  \label{optimum}
\Eeq
\Ethm

From now on, it is understood that the assumptions \HPstruttura\ and \gabri{those}
on the structure and on the initial data are satisfied
and that \pier{the} map~$\calS$, the cost functionals $\pier{\calF}$ and $\calJ$ 
and the control box $\Uad$ are defined \pier{in} \Defcontrol.
Thus, we do not remind anything of that in the statements given in the sequel.

Our next aim is to formulate \gabri{the} first order necessary optimality conditions.
As $\Uad$ is convex, the desired necessary condition for optimality~is
\Beq
  \< D\calJ(\uopt) , u-\uopt > \geq 0
  \quad \hbox{for every $u\in\Uad$}
  \label{precnopt}
\Eeq
provided that the derivative $D\calJ(\uopt)$ exists in the dual space $\calX^*$
at least in the G\^ateaux sense.
Then, the natural approach consists in proving that 
$\calS$ is \Frechet\ differentiable at $\uopt$
and applying the chain rule to $\calJ=\pier{\calF}\circ\calS$.
We can properly tackle this project under further assumptions 
on the nonlinearities~$\beta$, $\pi$ and~$g$.
Namely, we also suppose that
\Bsist
  & \hbox{$D(\beta)$ is an open interval and $\beta$ is a single-valued on $D(\beta)$}
  \label{hpDbetabis}
  \\
  & \hbox{$\beta$ and $\pi$ are $C^2$ functions
    and $g$ is a $C^1$ function}
  \label{hpnonlinbis}
\Esist
\Accorpa\HPstrutturabis hpDbetabis hpnonlinbis
\pier{and observe that, in particular,} $\betaz=\beta$.

We remark that both the regular potential \eqref{regpot} and the logarithmic potential \eqref{logpot}
satisfy the above assumptions on $\beta$ and~$\pi$.
Another possible choice of $\beta$ is given by
\Beq
  \beta(r) := 1 - \frac 1{\pier{r+1}}
  \quad \hbox{for $r>\pier{{}-1}$}
  \label{altrobeta}
\Eeq
and it corresponds to $\Beta$ defined by
\Beq
  \Beta(r) := r - \pier{\ln (r+1)}
  \quad \hbox{if $r>-1$}
  \aand
  \Beta(r) := +\infty
  \quad \hbox{otherwise}
  \label{altroBeta}
\Eeq
\pier{with $\Beta$ taking the minimum $0$ \gianni{at} $0$, as required by} assumption~\eqref{hpBeta}.
\pier{Such an operator $\beta$ yields an example of a different \bhv\ for negative and positive values, 
singular near $-1$ and with a somehow linear growth at $+\infty$.} 

Furthermore, we notice that the \pier{inclusion in} \eqref{seconda} becomes $\xi=\beta(\phi)$
and that $\beta$ and $\pi$ enter the problem through their sum, mainly.
Hence, we set for brevity
\Beq
  \gamma := \beta + \pi
  \label{defgamma}
\Eeq
and observe that $\gamma$ is a $C^2$ function on~$D(\beta)$.

\pier{Since} assumptions \HPstrutturabis\ force $\beta(r)$ to tend to $\pm\infty$
as $r$ tends to a finite end-point of~$D(\beta)$, if any,
we see that combining the further requirement \HPstrutturabis\
with the boundedness of $\phi$ and $\xi$ given by Theorem~\ref{Regularity}
immediately yields

\Bcor
\label{Bddaway}
Under all the assumptions of Theorem~\ref{Regularity},
suppose that \HPstrutturabis\ hold, in addition.
Then, the component $\phi$ of the solution $(\theta,\phi,\xi)$
also satisfies
\Beq
  \phimin \leq \phi \leq \phimax
  \quad \hbox{in $\overline Q$}
  \label{bddaway}
\Eeq
for some constants $\phimin\,,\phimax\in D(\beta)$
that depend only on $\Omega$, $T$, 
the structure \HPstruttura\ and \HPstrutturabis\ of the system,
the norms of the initial data associated to~\HPdati,
and the norms
$\norma u_\infty$, $\norma\thetaz_\infty$ and $\norma{\beta(\phiz)}_\infty$.
\Ecor

As we shall see in Section~\ref{FRECHET}, 
the computation of the \Frechet\ derivative of $\calS$ 
leads to the linearized problem that we describe at once
and that can be stated starting from a generic element $\ubar\in\calX$.
Let $\ubar\in\calX$ and $h\in\calX$ be given.
We set $(\thetabar,\phibar):=\calS(\ubar)$.
Then the linearized problem consists in finding $(\Theta,\Phi)$ 
satisfying 
\Bsist
  & \Theta \in \H1H \cap \L\infty\Vz \cap \L2\Hdue 
  \label{regTheta}
  \\
  & \Phi \in \pier{\H1H \cap \L\infty V \cap \L2 W}
  \label{regPhi}
\Esist
and solving the following problem
\Bsist
  & \dt\Theta - \Delta\Theta + \dt\Phi = mh
  & \quad \aeQ
  \label{linprima}
  \\
  & \dt\Phi - \Delta\Phi + \gamma'(\phibar) \, \Phi = \Theta
  & \quad \aeQ
  \label{linseconda}
  \\
  & \Theta(0) = 0
  \aand
  \Phi(0) = 0
  & \quad \aeO .
  \label{lincauchy}
\Esist
\Accorpa\Reglin regTheta regPhi
\Accorpa\Linpbl linprima lincauchy

\Bprop
\label{Existlin}
Let $\ubar\in\calX$ and $(\thetabar,\phibar)=\calS(\ubar)$.
Then, for every $h\in\calX$,
there exists a unique pair $(\Theta,\Phi)$
satisfying \Reglin\
and solving the linearized problem \Linpbl.
Moreover, the inequality
\Beq
  \norma{(\Theta,\Phi)}_\calY \leq C_5 \norma h_\calX
  \label{stimaFrechet}
\Eeq
holds true with a constant $C_5$ 
that depend only on $\Omega$, $T$, 
the structure \HPstruttura\ and \HPstrutturabis\ of the system,
the norms of the initial data associated to~\HPdati,
and the norms
$\norma\ubar_\infty$, $\norma\thetaz_\infty$ and $\norma{\beta(\phiz)}_\infty$.
In particular, the linear map $\calD:h\mapsto(\Theta,\Phi)$
is continuous from $\calX$ to~$\calY$.
\Eprop

Namely, we shall prove that the \Frechet\ derivative 
$D\calS(\ubar)\in\calL(\calX,\calY)$ 
actually exists and coincides with the map $\calD$ introduced in the last \pier{statement}.
This will be done in Section~\ref{FRECHET}.
Once this is established, we may use the chain rule with $\ubar:=\uopt$ 
to prove that the necessary condition \eqref{precnopt} for optimality takes the form
\Beq
  \intQ \bigl( g(\phiopt) - \chi \bigr) g'(\phiopt) \Phi
  + \kappa \intQ (\thetaopt - \thetaQ) \Theta
  \geq 0
  \quad \hbox{for any $u\in\Uad$},
  \qquad
  \label{cnopt}
\Eeq
where $(\thetaopt,\phiopt)=\calS(\uopt)$ and, for any given $u\in\Uad$, 
the pair $(\Theta,\Phi)$ is the solution to the linearized problem
corresponding to $h=u-\uopt$.

The final step then consists in eliminating the pair $(\Theta,\Phi)$ from~\eqref{cnopt}.
This will be done by introducing a pair $(p,q)$ 
that fulfills the regularity requirements
\Bsist
  && p \in \H1H \cap \L\infty\Vz \cap \L2{\Hx2}
  \label{regp}
  \\
  && q \in \H1H \cap \L\infty V \cap \L2W
  \label{regq}
\Esist
\Accorpa\Regadj regp regq
and solves the following adjoint system:
\Bsist
  & - \dt p - \Delta p - q
  = \kappa (\thetaopt-\thetaQ)
  & \quad \aeQ
  \label{primaadj}
  \\
  & -\dt q - \Delta q + \gamma'(\phiopt) \, q
  -  \dt p
  = \bigl( g(\phiopt) - \chi \bigr) g'(\phiopt)
  & \quad \aeQ
  \label{secondaadj}
  \\
  & p(T) = q(T) =0
  & \quad \aeO .
  \label{cauchyadj}
\Esist
\Accorpa\Pbladj primaadj cauchyadj
\pier{Here, let us recall \eqref{defspazi} and note that, as in previous cases (cf.~\accorpa{regtheta}{cauchy} and~\accorpa{regTheta}{lincauchy}), the Dirichlet boundary 
condition for $p$ is contained in \eqref{regp} whereas the
Neumann boundary condition for $q$ is in \eqref{regq}.}

\Bthm
\label{Existenceadj}
Let $\uopt$ and $(\thetaopt,\phiopt)=\calS(\uopt)$
be an optimal control and the corresponding state.
Then the adjoint problem \Pbladj\ has a unique solution $(p,q)$
satisfying the regularity conditions \Regadj.
\Ethm

\pier{Our last result establishes optimality conditions.}

\Bthm
\label{CNoptadj}
Let $\uopt$ be an optimal control.
Moreover, let $(\thetaopt,\phiopt)=\calS(\uopt)$ and $(p,q)$
be the associate state and the unique solution to the adjoint problem~\Pbladj\
given by Theorem~\ref{Existenceadj}.
Then we~have\pier{%
\begin{align}
   m(x) \, p(t,x) \, \bigl( \gabri{u}-\uopt(t,x) \bigr) \geq 0 
   \quad \hbox{for every $\gabri{u}\in[\umin(t,x),\umax(t,x)],$} 
   \quad \non
   \\
   \quad \aaQ .   
   \label{cnoptadj}
\end{align}
}%
In particular, $mp=0$ in the subset of $Q$ where $\umin<\uopt<\umax$.
\Ethm

\pier{A straightforward consequence of Theorem~\ref{CNoptadj} is here stated.}

\Bcor
\gabri{Under the conditions of Theorem~\ref{CNoptadj}, the
optimal control $\uopt $ reads}
\[\gabri{
u^{\ast }=\left\{
\begin{array}{ll}
u_{\min } &\mbox{a.e. on the set } \ \{(t,x) \, : \ p(t,x)>0\mbox{
and }m(x)>0\} \\[0.1cm]
u_{\max } &\mbox{a.e. on the set } \ \{(t,x)\, : \ p(t,x)<0\mbox{
and }m(x)>0\} \\[0.1cm]
\mbox{undetermined}\quad &\mbox{elsewhere.}%
\end{array}%
\right. }
\]
\Ecor

In the remainder of the paper, we often owe to the \Holder\ inequality
and to the elementary Young inequalities
\Bsist
  && ab \leq \alpha \, a^{1/\alpha} + (1-\alpha) \, b^{1/(1-\alpha)}
  \aand
  ab \leq \delta a^2 + \frac 1 {4\delta} \, b^2
  \non
  \\
  && \quad \hbox{for every $a,b\geq 0$, \ $\alpha\in(0,1)$ \ and \ $\delta>0$}
  \label{young}
\Esist
in performing our a priori estimates.
To this regard, in order to avoid a boring notation,
we use the following general rule to denote constants.
The small-case symbol $c$ stands for different constants which depend only
on~$\Omega$, on the final time~$T$, the shape of the nonlinearities
and on the constants and the norms of
the functions involved in the assumptions of our statements.
A~small-case $c$ with a subscript like $c_\delta$
indicates that the constant might depend on the parameter~$\delta$, in addition.
Hence, the meaning of $c$ and $c_\delta$ might
change from line to line and even in the same chain of equalities or inequalities.
On the contrary, different symbols (e.g., capital letters)
stand for precise constants which we can refer~to.


\section{The state system}
\label{STATE}
\setcounter{equation}{0}

This section is devoted to \betti{the proofs of} Theorems~\ref{Wellposedness} and~\ref{Regularity}.
As far as the former is concerned,
we notice that the initial--boundary value problem under study 
is a quite standard phase field system
and that a number of results on it can be found in the literature
(see, e.g., \pier{\cite{BrokSpr, DKS, Sch}}, and references therein).
Nevertheless, we prefer to sketch the basic a~priori estimates
that correspond to the regularity \Regsoluz\ of the solution
and to estimate~\eqref{stimasoluz},
for the reader convenience.
A~complete existence proof can be obtained by \regulariz ing the problem,
performing the same estimates on the corresponding solution,
and passing to the limit through compactness results.
We also give a short proof of~\eqref{contdep}
(whence \pier{uniqueness follows} as a consequence)
and conclude the discussion on Theorem~\ref{Wellposedness}.
 
As said, we derive just formal a priori estimates.
\pier{We multiply \eqref{prima} by $\theta$; then we add $\phi$ to both sides
of \eqref{seconda} and test by $\dt\phi$;
finally, we sum up and integrate over~$Q_t$ with $t\in(0,T)$.
As the terms involving the product $\theta\,\dt\phi$ cancel out, 
by exploiting a standard chain rule for subdifferentials (see, 
e.g., \cite[Lemme~3.3, p.~73]{Brezis}) we obtain
\begin{align}
  & \frac 12 \iO |\theta(t)|^2  
  + \intQt |\nabla\theta|^2
  +  \intQt |\dt\phi|^2
  + \frac 12 \norma{\phi (t)}_V^2
  + \frac  12 \iO \Beta(\phi(t))
  \non
  \\
  & = \frac 12 \iO |\thetaz|^2
  + \frac  12 \norma{\phiz}_V^2
   + \frac  12 \iO \Beta(\phiz)
     + \intQt m u \, \theta
  +  \intQt (\phi - \pi (\phi)) \, \dt\phi .
  \label{p5}
\end{align}
The last integral on the \lhs\ is nonnegative thanks to \eqref{hpBeta}
and the first three terms on the \rhs\ are under control, due to \accorpa{hpz}{hpzbis}.
Since (cf.~\accorpa{hpPi}{defbetapi}) 
$ |\phi- \pi (\phi)| \leq c(|\phi| + 1)  $ and \eqref{hplm} holds, 
the last two terms on the \rhs\ of \eqref{p5} can be easily dealt with 
by the Young inequality and the Gronwall lemma. Then, we deduce the estimate 
 \Beq
  \norma\theta_{\L\infty H\cap\L2 V}
  + \norma\phi_{\H1H\cap\L\infty V}
  \leq c . 
  \label{p6}
\Eeq 
Since $\dt\phi$ is by now bounded in $L^2(Q)$, we can 
test \eqref{prima} by $\dt\theta$ in order to infer that 
\Beq
  \intQt |\dt\theta|^2
  + \frac 12 \iO |\nabla\theta(t)|^2
  = \frac 12 \iO |\nabla\thetaz|^2
  + \intQt (m u - \dt\phi) \, \dt\theta .
  \non
  \Eeq
Thus, \eqref{hpz} and the Young inequality enable us to recover
\Beq
  \norma\theta_{\H1H\cap\L\infty V}
  \leq c 
  \label{p7}
\Eeq
as well. At this point, owing to \eqref{p6}--\eqref{p7}, 
$\Delta\theta$ and $-\Delta\phi+\xi$ are bounded 
in $\L2H$, as one clearly sees from equations \accorpa{prima}{seconda}.
Hence, a~standard monotonicity argument (test some regularization of \eqref{seconda} by the analogue of $\xi = \beta(\phi)$) yields that both $\Delta\phi$ and $\xi$ are bounded in $\L2 H$. Then, elliptic regularity allows us to derive the complete estimate~\eqref{stimasoluz}.
}

Let us pass to~\eqref{contdep}.
We first integrate \eqref{prima} with respect to time
and get the equation
\Beq
  \theta - \Delta(1*\theta) + \phi
  = \thetaz + \phiz + m (1*u) .
  \label{intprima}
\Eeq
Now, we fix $u_i\in\LQ2$, $i=1,2$,
and consider two corresponding solutions $(\theta_i,\phi_i,\xi_i)$ \betti{with the same initial data}.
We write \eqref{intprima} for both of them 
and multiply the difference by $\theta:=\theta_1-\theta_2$.
At the same time, we write \eqref{seconda} for both solution 
and multiply the difference by $\phi$, where $\phi:=\phi_1-\phi_2$.
Then, we add the equalities we obtain to each other and integrate over~$Q_t$.
The terms involving the product $\phi\theta$ cancel out.
Hence, by also setting $u:=u_1-u_2$ and $\xi:=\xi_1-\xi_2$ for brevity, we have
\Bsist
  && \intQt |\theta|^2
  + \frac 12 \iO |\nabla(1*\theta)(t)|^2
  + \frac  12 \iO |\phi(t)|^2
  +  \intQt |\nabla\phi|^2
  +  \intQt \xi \phi
  \non
  \\
  && = \intQt m (1*u) \, \theta
  - \intQt \bigl( \pi(\phi_1) - \pi(\phi_2) \bigr) \, \phi
  \non
  \\
  && \leq c \norma{1*u}_{\LQ2}^2 + \frac 12 \intQt |\theta|^2
  + c \intQt |\phi|^2 
  \non
\Esist 
where we used the boundedness of $m$ and the \Lip\ continuity of~$\pi$
(see \eqref{hplm} and \accorpa{hpPi}{defbetapi} once more).
As the last integral on the \lhs\ is nonnegative since $\beta$ is monotone,
we obtain the desired estimate \eqref{contdep} 
just by rearranging and applying the Gronwall lemma.\qed

\medskip

Now, we prove Theorem~\ref{Regularity} \pier{using the same strategy 
of a formal argumentation. First, we consider the equation obtained by differentiating \eqref{seconda} with respect to time and test it by $\dt \phi$. Then, we have  
\Bsist
  &&
  \frac  12 \iO |\dt\phi(t)|^2
  +  \intQt |\nabla\dt\phi|^2
  +  \intQt \beta'(\phi) |\dt\phi|^2
  \non
  \\
  && \leq 
  \frac  12 \iO |\dt\phi(0)|^2
  + \intQt (\dt\theta - \pi'(\phi) \dt\phi) \dt\phi . 
  \non
\Esist
The monotonicity of $\beta$ implies that the last term on the \lhs\ is nonnegative;
on the \rhs , the last integral is already bounded thanks to \accorpa{p6}{p7} 
and to the boundedness of  $\pi'$ (see \accorpa{hpPi}{defbetapi}).
Thus, just the norm of $\dt\phi(0)$ in $\Lx2$ should be estimated,
and this can be performed by recovering $\dt\phi(0)$
from equation \eqref{seconda} and then exploiting~\HPdati\ as well as \eqref{p1}.
Consequently, we obtain 
\Beq
 \norma\phi_{\W{1,\infty}H\cap\H1V}
  \leq c 
  \non
\Eeq
and, in addition, the boundedness of $-\Delta\phi+\xi$ in $\L\infty H$. 
Now, it is straightforward to infer that
the two separate terms $\Delta\phi$ and $\xi$ are both bounded in $\L\infty H$.
Then, the properties \accorpa{p2}{p4} follow; moreover, they} imply that 
$\phi$ is bounded in $\C0{\Cx0}=C^0(\overline Q)$
since $W$ is complactly embedded in~$\Cx0$
(see, e.g., \cite[Sect.~8, Cor.~4]{Simon}).
This proves~$i)$.
For the second statement~$ii)$, 
we observe that $\theta$~turns out to be bounded whenever its initial value is bounded.
Indeed, \eqref{seconda} can be written in the form
\Beq
  \dt\theta - \Delta\theta = mu - \dt\phi \in \L\infty H
  \non
\Eeq
whence it suffices to apply, e.g.,
\cite[\pier{Thm.}~7.1, p.~181]{LSU} with $r=\infty$ and $q=2$.
Finally, we prove~$iii)$ by writing \eqref{prima} in the form
\Beq
  \dt\phi - \Delta\phi + \xi = f:= \theta - \pi(\phi) 
  \aand
  \xi \in \beta(\phi)
  \quad \aeQ
  \label{secondabis}
\Eeq
and \gianni{observing} that $f$ is bounded in~$\LQ\infty$
on account of the result $ii)$ just proved.
Now, we approximate $\phi$ by the solution $\phieps$
to the initial--boundary value problem obtained 
by keeping the same initial and boundary conditions and replacing \eqref{secondabis}~by
\Beq
  \dt\phieps - \Delta\phieps + \xieps = f:= \theta - \pi(\phi) 
  \aand
  \xieps := \betaeps(\phieps)
  \quad \aeQ
  \label{secondaeps}
\Eeq
where $\betaeps$ is the Yosida \regulariz ation of~$\beta$ at level~$\eps>0$.
Indeed, a~standard argument shows that $\phieps$ converges to $\phi$ 
in the proper topology as $\eps$ tends to zero,
so that $iii)$~immediately follows
whenever we prove that $\xieps$ is bounded in $\LQ\infty$ uniformly with respect to~$\eps$.
To this end, by extending the sign function by $\sign(0)=0$,
we notice that $\sign\betaeps(r)=\sign r$ for every $r\in\erre$ 
since $\beta(0)\ni0$ (see~\eqref{hpBeta})
and~set
\Beq
  \Betap(r) := \int_0^r |\betaeps(s)|^{p-1} \sign s \, ds
  \quad \hbox{for $r\in\erre$ and $p>1$}.
  \non
\Eeq
We obtain a nonnegative function.
Then, we multiply \eqref{secondaeps} by $|\xieps|^{p-1}\sign\xieps$,
where $p>2$ is arbitrary, and integrate over~$Q_t$.
We have
\Bsist
  && \iO \Betap(\phieps(t))
  + (p-1) \intQt |\xieps|^{p-2} \betaeps'(\phieps) |\nabla\phieps|^2 
  + \intQt |\xieps|^p
  \non
  \\
  && = \iO \Betap(\phiz)
  + \intQt f \, |\xieps|^{p-1} \sign\xieps .
  \non
\Esist
By noting that the first \pier{two} terms on the \lhs\ are nonnegative
and owing to the Young inequality, we deduce that
\Beq
  \intQt |\xieps|^p
  \leq \iO \Betap(\phiz)
  + \intQt |f| \, |\xieps|^{p-1} 
  \leq \iO \Betap(\phiz)
  + \frac 1p \intQt |f|^p + \frac 1{p'} \intQt |\xieps|^p .
  \non
\Eeq
By rearranging, we obtain
\Beq
  \intQt |\xieps|^p
  \leq p \iO \Betap(\phiz)
  + \intQt |f|^p
  \non
\Eeq
whence also (since $(a+b)^\alpha\leq a^\alpha+b^\alpha$ for every $a,b\geq0$ and $\alpha\in(0,1)$)
\Beq
  \norma\xieps_{\LQ p}
  \leq p^{1/p} \Bigl( \iO \Betap(\phiz) \Bigr)^{1/p}
  + (|\Omega|T)^{1/p} \norma f_\infty \,.
  \non
\Eeq
By letting $p$ tend to infinity, we conclude that
\Beq
  \norma\xieps_\infty \leq \cz + \norma f_\infty
  \quad \hbox{provided that} \quad
  \iO \Betap(\phiz) \leq \cz^p
  \non
\Eeq
and we just have to show that such a finite $\cz$ exists.
To this aim, we notice that $r$, $\betaeps(r)$ and $\betaz(r)$
have the same sign for every $r\in\erre$\pier{; on the other hand, 
\eqref{propYosida} holds
and even} $\betaz$ is monotone.
Hence, we have 
\Beq
  \Betap(\phiz)
  = \Bigl| \int_0^{\phiz} |\betaeps(s)|^{p-1} \, ds \Bigr|
  \leq |\phiz| \, |\betaz(\phiz)|^{p-1} \pier{\quad \aeO.}
  \non
\Eeq
As both $\phiz$ and $\betaz(\phiz)$
are bounded, the former since $\phiz\in W\subset\Lx\infty$
and the latter by assumption, we deduce~that
\Beq
  \iO \Betap(\phiz)
  \leq |\Omega| \, \norma\phiz_\infty \, \norma{\betaz(\phiz)}_\infty^{p-1}
  \leq \cz^p
  \non
\Eeq
with an obvious choice of~$\cz$, and the proof is complete.\qed


\section{Existence of an optimal control}
\label{OPTIMUM}
\setcounter{equation}{0}

\gabri{The following section is devoted to the proof of  Theorem~\ref{Optimum}. We use} the direct method,
observing \betti{first} that $\Uad$ is nonempty.
\betti{Then, we} let $\graffe{\un}$ be a minimizing sequence for the optimization problem
and, for any~$n$, \betti{we} take the corresponding solution $(\phin,\thetan,\xin)$ to problem~\Pbl.
Then, $\graffe{\un}$ is bounded in $\Lx\infty$ and estimate \eqref{stimasoluz} holds for $(\phin,\thetan,\xin)$.
Therefore, we have for a subsequence
\Bsist
  & \un \to u
  & \quad \hbox{weakly star in $\Lx\infty$}
  \non
  \\
  & \thetan \to \theta
  & \quad \hbox{weakly star in $\H1H \cap \L\infty\Vz \cap \L2\Hdue$}
  \non
  \\
  & \phin \to \phi
  & \quad \hbox{weakly star in $\W{1,\infty}H \cap \H1V \cap \L\infty W$}
  \non
  \\
  & \xin \to \xi
  & \quad \hbox{weakly star in $\L\infty H$} .
  \non
\Esist
Then, $u\in\Uad$ since $\Uad$ is closed in~$\calX$,
the initial conditions for $\theta$ and $\phi$ are satisfied,
and we can easily conclude by standard argument.
Very shortly, $\graffe{\phin}$~converges strongly, e.g., in~$\pier{L^2(Q)}$ and \aeQ\ (for~a subsequence)
by~the Aubin-Lions compactness lemma (see, e.g., \cite[Thm.~5.1, p.~58]{Lions}),
whence $\pi(\phin)$ converges to $\pi(\phi)$ is the same topology
and $\xi\in\beta(\phi)$ (see, e.g., \cite[Lemma~1.3, p.~42]{Barbu}).
Thus, $(\theta,\phi,\psi)$ satisfies the variational formulation in the integral form of problem~\Pbl.
On the other hand, $\pier{\calF}(\thetan,\phin)$ converges both to the infimum of $\calJ$
and to $\pier{\calF}(\theta,\phi)$,
since $g(\phin)$ converges to $g(\phi)$ \aeQ\ and \pier{it} is bounded in $\LQ\infty$ 
(see~\eqref{hpg}).
Therefore, $u$~is an optimal control.\qed


\section{The control-to-state mapping}
\label{FRECHET}
\setcounter{equation}{0}

As sketched in Section~\ref{STATEMENT}, 
the main point is the \Frechet\ differentiability of the control-to-state mapping~$\calS$.
This involves the \lineariz ed problem \Linpbl\pier{,
whose} well-posedness is stated in Proposition~\ref{Existlin}.
Thus, we first prove such a result.

As one can easily see by going through the proof of estimates
\eqref{stimasoluz} and~\eqref{contdep} given in Section~\ref{STATE},
what is stated in Theorem~\ref{Wellposedness} can be extended to the problem
obtained by replacing equation \eqref{seconda} of \Pbl\ by the more general~one
\Beq  
  \dt\phi - \Delta\phi + \xi + \alpha\pi(\phi) = \theta
  \aand
  \xi \in \beta(\phi)
  \non
\Eeq
where $\alpha\in\pier{\LQ\infty}$ is prescribed.
Therefore, Proposition~\ref{Existlin} follows as a trivial particular case.
Namely, one just chooses $\beta=0$, $\pi(r)=r$ and $\alpha=\gamma'(\phibar)$
where $\gamma$ is defined by \eqref{defgamma} by starting from the original $\beta$ and~$\pi$.
Indeed, $\gamma'(\phibar)$~is bounded thanks to Corollary~\ref{Bddaway}.
In fact, estimate \eqref{stimaFrechet} holds more generally with $\norma h_{\LQ2}$ on the \rhs.\qed

Here is the main result of this section.

\Bthm
\label{Fdiff}
Let $\ubar\in\calX$ and $(\thetabar,\phibar)=\calS(\ubar)$.
Then, $\calS$~is \Frechet\ differentiable at~$(\thetabar,\phibar)$
and the \Frechet\ derivative $[D\calS](\thetabar,\phibar)$ 
precisely is the map $\calD\in\calL(\calX,\calY)$ defined
in the statement of Proposition~\ref{Existlin}.
\Ethm

\Bdim
We fix $\ubar\in\calX$ and the corresponding state $(\thetabar,\phibar)$
and, for $h\in\calX$ with $\norma h_\calX\leq\gabri{\Lambda}$\gabri{, for some positive constant $\Lambda$,} we~set
\Beq
  (\thetah,\phih) := \calS(\pier{\ubar}+h)
  \aand
  (\zetah,\etah) := (\thetah-\thetabar-\Theta,\phih-\phibar-\Phi)
  \non
\Eeq
where $(\Theta,\Phi)$ is the solution to the linearized problem corresponding to~$h$.
We have to prove that
$\norma{(\zetah,\etah)}_\calY/\norma h_\calX$ tends to zero as $\norma h_\calX$ tends to zero.
More precisely, we show that
\Beq
  \norma{(\zetah,\etah)}_\calY
  \leq c \norma h_{\LQ2}^2
  \label{tesiFrechet}
\Eeq
for some constant~$c$,
and this is even stronger than necessary.
First of all, we fix one fact.
As both $\norma\ubar_\infty$ and $\norma{\ubar+h}_\infty$
are bounded by $\norma\ubar_\infty+\gabri{\Lambda}$,
we can apply Corollary~\ref{Bddaway}
and find constants $\phimin,\phimax\in D(\beta)$ such that
\Beq
  \phimin \leq \phibar \leq \phimax
  \aand
  \phimin \leq \phih \leq \phimax
  \quad \aeQ .
  \label{perstimaTaylor}
\Eeq
Now, let us prove~\eqref{tesiFrechet} by writing the problem solved by~$(\zetah,\etah)$.
We clearly have
\Bsist
  & \dt\zetah - \Delta\zetah +  \dt\etah = 0
  & \quad \aeQ
  \label{primah}
  \\
  & \dt\etah - \Delta\etah + \gamma(\phih) - \gamma(\phibar) - \gamma'(\phibar) \, \Phi = \zetah
  & \quad \aeQ .
  \label{secondah}
\Esist
Moreover, both $\zetah$ and $\etah$ satisfy homogeneous initial and boundary conditions
(of~Dirichlet and Neumann type, respectively).
Now, we integrate \eqref{primah} with respect to time and obtain
\Beq
  \zetah - \Delta (1*\gabri{\zetah}) +  \etah = 0 .
  \label{intprimah}
\Eeq
At this point, we multiply \eqref{intprimah} and \eqref{secondah} 
by $\zetah$ and $\etah$, respectively, 
add the resulting equalities to each other and integrate over~$Q_t$.
The terms involving the product $\zetah\etah$ cancel out and we have
\Bsist
  && \intQt |\zetah|^2
  + \frac 12 \iO |(1*\nabla\zetah)(t)|^2
  + \frac 1 2 \iO |\etah(t)|^2
  +  \intQt |\nabla\etah|^2
  \non
  \\
  && = -  \intQt \bigl( \gamma(\phih) - \gamma(\phibar) - \gamma'(\phibar) \, \Phi \bigr) \etah \,.
  \label{dastimare}
\Esist
Now, \aaQ, we write the Taylor expansion of $\gamma$ around $\phibar(t,x)$.
Some function $\tilde\phi_h$ exists such that
\Bsist
  && \gamma(\phih)
  = \gamma(\phibar)
  + \gamma'(\phibar) (\phih-\phibar)
  + \frac 12 \, \gamma''(\tilde\phi_h) (\phih-\phibar)^2
  \quad \aeQ
  \non
  \\
  && \min\{\phih,\phibar\} \leq \tilde\phi_h \leq \max\{\phih,\phibar\}
  \quad \aeQ .
  \non
\Esist
Then, $\phimin\leq\tilde\phi_h\leq\phimax$ by \eqref{perstimaTaylor}.
It follows that $\gamma''(\tilde\phi_h)$ is bounded
since $D(\beta)$ is an open interval and $\gamma''$ is continuous.
As the same is true for~$\gamma'(\phibar)$, we can estimate the \rhs\ of \eqref{dastimare}
by accounting for the Young and \Holder\ inequalities with any $\delta\in(0,1)$ as follows
\Bsist
  && -  \intQt \bigl( \gamma(\phih) - \gamma(\phibar) - \gamma'(\phibar) \, \Phi \bigr) \etah
  \non
  \\
  && = -  \intQt \bigl( \gamma'(\phibar) \etah + \frac 12 \, \gamma''(\tilde\phi_h) (\phih-\phibar)^2 \bigr) \etah
  \non
  \\
  && \leq c \intQt |\etah|^2
  + c \intQt |\phih-\phibar|^2 |\etah|
  \non
  \\
  && \leq c \intQt |\etah|^2
  + c \iot \norma{(\phih-\phibar)(s)}_4 \, \norma{(\phih-\phibar)(s)}_2 \, \norma{\etah(s)}_4 \, ds
  \non
  \\
  && \leq c \intQt |\etah|^2
  + \delta \iot \norma{\etah(s)}_4^{\betti{2}} \, ds
  + c_\delta \iot \norma{(\phih-\phibar)(s)}_4^2 \, \norma{(\phih-\phibar)(s)}_2^2 \, ds .
  \non
\Esist
Now, we recall the Sobolev inequality $\norma v_4\leq C_\Omega\normaV v$ for every $v\in V$, 
where $C_\Omega$ depends only on~$\Omega$,
and that estimate \eqref{contdep} holds for the pair of controls $\ubar+h$ and $\ubar$
and for the corresponding states 
$(\thetah,\phih)$ and $(\thetabar,\phibar)$.
Therefore, we can continue and obtain
\Bsist
  && -  \intQt \bigl( \gamma(\phih) - \gamma(\phibar) - \gamma'(\phibar) \, \Phi \bigr) \etah
  \non
  \\
  && \leq c \intQt |\etah|^2
  + \delta \, C_\Omega \intQt \bigl( |\etah|^2 + |\nabla\etah|^2 \bigr)
  + c_\delta \norma{\phih-\phibar}_{\L2V}^2 \, \norma{\phih-\phibar}_{\L\infty H}^2
  \non
  \\
  && \leq c \intQt |\etah|^2
  + \delta \, C_\Omega \intQt |\nabla\etah|^2
  + c_\delta \norma h_{\LQ2}^4 \,.
  \non
\Esist
At this point, we choose $\delta$ small enough, rearrange and apply the Gronwall lemma.
This yields~\eqref{tesiFrechet}.
\Edim


\section{Necessary optimality conditions}
\label{OPTIMALITY}
\setcounter{equation}{0}

In this section, we derive the optimality condition~\eqref{cnoptadj}
stated in Theorem~\ref{CNoptadj}.
We start from~\eqref{precnopt} and first prove~\eqref{cnopt}.

\Bprop
\label{CNopt}
Let $\uopt$ be an optimal control and $(\phiopt,\thetaopt):=\calS(\uopt)$.
Then, \eqref{cnopt} holds.
\Eprop

\Bdim
This is \pier{essentially} due to the chain rule for \Frechet\ derivatives, 
as already said in Section~\ref{STATEMENT},
and we just provide some detail.
We notice that $g$ and $g'$ are computed 
only at the values of $\phiopt$ in~\eqref{cnopt}
and we can fix $\phimin,\phimax\in D(\beta)$
in order that \eqref{bddaway} holds for~$\phiopt$
and \pier{modify} $g$ outside of $[\phimin,\phimax]$ 
without changing anything else both in the problem and in the formula we want to prove.
Hence, we can assume even $g'$ to be bounded
so that the functional
\Beq
  \phi \mapsto \frac 12 \intQ \bigl( g(\phi) - \chi)^2
  \non
\Eeq
is well-defined and \Frechet\ differentiable in the whole of $\LQ2$.


It follows that $\pier{\calF}$ is \Frechet\ differentiable in $\calZ:=\LQ2\times\LQ2$
and that its \Frechet\ derivative $[D\pier{\calF}](\thetabar,\phibar)$ 
at any point $(\thetabar,\phibar)\in\calZ$ acts as follows
\Beq
  [D\pier{\calF}](\thetabar,\phibar):
  (h_1,h_2)\in\calZ \mapsto
  \intQ \bigl( g(\phibar) - \chi \bigr) g'(\phibar) h_1
  + \kappa \intQ (\thetabar - \thetaQ) h_2 \,.
  \non
\Eeq
Therefore, Theorem~\ref{Fdiff} and the chain rule ensure that
$\calJ$ is \Frechet\ differentiable at $\uopt$
and that its \Frechet\ derivative $[D\calJ](\uopt)$ and any optimal control $\uopt$
acts as follows
\Beq
  [D\calJ](\uopt):
  h \in \calX \mapsto
  \intQ \bigl( g(\phiopt) - \chi \bigr) g'(\phiopt) \Phi
  + \kappa \intQ (\thetaopt - \thetaQ) \Theta
  \non
\Eeq
where $(\Theta,\Phi)$ is the solution to the linearized problem corresponding to~$h$.
Therefore, \eqref{cnopt} immediately follows from~\eqref{precnopt}.
\Edim

The next step is the proof of Theorem~\ref{Existenceadj}.
For convenience, we consider the equivalent forward problem 
in the unknown $(\ptilde,\qtilde)$
given by $(\ptilde,\qtilde)(t):=(p,q)(T-t)$.
However, to simplify notations, we write $p$ and $q$ instead 
of $\ptilde$ and~$\qtilde$ in the sequel.
\pier{Thus, we write the homogeneous initial--boundary value problem 
\Bsist
  & \dt p - \Delta p - q = f_1
  & \quad \aeQ
  \label{primaf}
  \\
  & \dt q - \Delta q + \alpha q +  \dt p = f_2 
  & \quad \aeQ
  \label{secondaf}
  \\
  & p = 0 
  \aand
  \dn q = 0
  & \quad \aeS
  \label{Ibcf}
  \\
  & p(0) = 0
  \aand
  q(0) = 0
  & \quad \aeO 
  \label{Icauchyf}
\Esist
}%
with an obvious choice of $f_1,f_2\in\LQ2$ and $\alpha\in\LQ\infty$.
In order to prove uniqueness, we replace $f_1$ and $f_2$ by~$0$.
We multiply the above equations by $\dt p$ and~$q$, respectively,
add the equalities we get to each other
and observe that the terms involving the product $q\dt p$ cancel out.
Then, we integrate over $Q_t$ and rearrange.
We obtain
\Beq
   \intQt |\dt p|^2
  + \frac  12 \iO |\nabla p(t)|^2
  + \frac 12 \iO |q(t)|^2
  + \intQt |\nabla q|^2
  = - \intQt \alpha |q|^2 .
  \label{basicadj}
\Eeq
As $\alpha$ is bounded, we can apply the Gronwall lemma 
and conclude that $p=0$ and $q=0$.
As far as existence is concerned,
we start deriving the basic formal estimates.
As before, we multiply \eqref{primaf} by $\dt p$ and \eqref{secondaf} by~$q$
and perform the same calculation.
We obtain an inequality like \eqref{basicadj}
with a different \rhs, namely
\Beq
  - \intQt \alpha |q|^2
  +  \intQt f_1 \dt p
  + \intQt f_2 q \,.
  \non
\Eeq
By owing to the \Holder\ and Young inequalities,
we infer that the above expression is bounded from above~by
\Beq
  c \intQt |q|^2
  + \frac 1 2 \intQt |\dt p|^2
  + c \intQt \bigl( |f_1|^2 + |f_2|^2 \bigr).
  \non
\Eeq
Hence\pier{, we have that}
\Beq
  \norma p_{\H1H\cap\L\infty V}
  + \norma q_{\L\infty H\cap\L2V}
  \leq c
  \label{formala}
\Eeq
and the estimate
\Beq
  \norma{\dt q}_{\L2H}
  \leq c
  \label{formalb}
\Eeq
immediately follows as a consequence by multiplying \eqref{secondaf} by~$\dt q$.
Therefore, it is clear how to give a rigorous proof
based on a Faedo--Galerkin scheme,
which provides a sequence $\graffe{(p_n,q_n)}$ of approximating solutions
obtained by solving just linear systems of ordinary differential equations.
Namely, by performing the above estimates 
on $(p_n,q_n)$ exactly in the same way
and using standard compactness results,
one finds a weak limit $(p,q)$ in the topologies associated to \accorpa{formala}{formalb}
and it is immediately \pier{clear} that $(p,q)$ is a variational solution 
of the problem we want to solve.
Then, the complete regularity \Regadj\ and the fact that $(p,q)$
solves the problem in its strong form follow from the general theory.
So, Theorem~\ref{Existenceadj} actually holds.\qed

At this point, we are ready to prove Theorem~\ref{CNoptadj} on optimality,
i.e., the necessary condition \eqref{cnoptadj} for
$\uopt$ to be an optimal control in terms of the solution $(p,q)$
of the adjoint problem~\Pbladj.
So, we fix an arbitrary $u\in\Uad$
and write the variational formulations of both the \lineariz ed problem 
\pier{(}corresponding to $h=u-\uopt$\pier{)} and the adjoint problem.


The equations become
\Bsist
  && \intQ \dt\Theta \, v + \intQ \nabla\Theta \cdot \nabla v +  \intQ \dt\Phi \, v = \intQ m(u-\uopt) \, v
  \label{lvprima}
  \\
  && \intQ \dt\Phi \, v + \intQ \nabla\Phi \cdot \nabla v + \intQ \gamma'(\phiopt) \, \Phi \, v = \intQ \Theta \, v
  \label{lvseconda}
  \\
  && - \intQ \dt p \, v + \intQ \nabla p \cdot \nabla v - \intQ q v
  = \kappa \intQ (\thetaopt-\thetaQ) v
  \label{avprima}
  \\
  && - \intQ \dt q \, v + \intQ \nabla q \cdot \nabla v + \intQ \gamma'(\phiopt) \, q \, v
  -  \intQ \dt p \, v
  \non
  \\
  && = \intQ \bigl( g(\phiopt) - \chi \bigr) g'(\phiopt) v
  \label{avseconda}
\Esist
where \eqref{lvprima} and \eqref{avprima} have to hold for every $v\in\L2\Vz$,
while \eqref{lvseconda} and~\eqref{avseconda} are required for every $v\in\L2V$.
In particular \eqref{lvseconda} and~\eqref{avseconda} also contain the homogeneous Neumann conditions for $\Phi$ and~$q$.
Moreover, $\Theta$ and $p$ have to satisfy the homogeneous Dirichlet boundary conditions, in addition.
Finally, the functions at hand satisfy the homogeneous initial or final conditions 
as specified in \eqref{lincauchy} and~\eqref{cauchyadj}.
We choose $v=p$, $v=q$, $v=-\Theta$ and $v=-\Phi$ in \accorpa{lvprima}{avseconda}, respectively,
and add all the equality we obtain to each other.
The most part of the terms cancel out and we obtain
\Beq
  \intQ \dt \bigl( \Theta p + \Phi q +  \Phi p \bigr)
  = \intQ \bigl\{ m(u-\uopt) p - \kappa (\thetaopt-\thetaQ) \Theta - \bigl( g(\phiopt) - \chi \bigr) g'(\phiopt) \Phi \bigr\}\betti{\,.}
  \non
\Eeq
Due to the initial and final conditions, the \lhs\ vanishes and we deduce that
\Beq
  \intQ \bigl\{ \kappa (\thetaopt-\thetaQ) \Theta + \bigl( g(\phiopt) - \chi \bigr) g'(\phiopt) \Phi \bigr\}
  = \intQ m(u-\uopt) p \,.
  \non
\Eeq
As the \lhs\ is $\geq0$ by \eqref{cnopt}, it follows that the same is true for the \rhs.
As $u\in\Uad$ is arbitrary, this implies the \pier{pointwise inequality
 \eqref{cnoptadj}} and the proof of Theorem~\ref{CNoptadj} is complete.\qed


\pier{%
\section*{Acknowledgements}
This research activity has been performed in the framework of an
Italian-Romanian  \gabri{three-year project on ``Nonlinear partial differential equations (PDE) with applications in modeling cell growth, chemotaxis and phase transition'' financed by the Italian CNR and the Romanian Academy.} Moreover, 
the financial support of the FP7-IDEAS-ERC-StG \#256872
(EntroPhase) is gratefully acknowledged by the authors. The present paper 
also benefits from the support of the MIUR-PRIN Grant 2010A2TFX2 ``Calculus of Variations'' for PC and GG, \gabri{and} the GNAMPA (Gruppo Nazionale per l'Analisi Matematica, la Probabilit\`a e le loro Applicazioni) of INdAM (Istituto Nazionale di Alta Matematica) for PC, GG and ER.
}%


\vspace{3truemm}


\Begin{thebibliography}{10}

\bibitem{Barbu}
V. Barbu,
``Nonlinear semigroups and differential equations in Banach spaces'',
Noord-hoff, 
Leyden, 
1976.

\pier{
\bibitem{BBCG}
V. Barbu, M.L. Bernardi, P. Colli, G. Gilardi,
Optimal control problems of phase relaxation models,
J. Optim. Theory Appl. {\bf 109} (2001), 557-585.}

\pier{\bibitem{BCF} 
J.L. Boldrini, B.M.C. Caretta, E. Fern{\'a}ndez-Cara,
Some optimal control problems for a two-phase field model of solidification,
Rev. Mat. Complut. {\bf 23} (2010), 49-75.}

\bibitem{Brezis}
H. Brezis,
``Op\'erateurs maximaux monotones et semi-groupes de contractions
dans les espaces de Hilbert'',
North-Holland Math. Stud.
{\bf 5},
North-Holland,
Amsterdam,
1973.

\betti{
\bibitem{BrokSpr} 
\pier{M. Brokate, J. Sprekels,}
``Hysteresis and Phase Transitions'',
Springer, New York, 1996.
}

\betti{
\bibitem{Cag}
G.~Caginalp,
An analysis of a phase field model of a free boundary,
Arch. Rational Mech. Anal \pier{{\bf 92} (1986), 205-245}.
} 

\betti{
\bibitem{CGPS}
\pier{P. Colli, , G.~Gilardi,} P. Podio-Guidugli, J. Sprekels, 
Distributed optimal control of a nonstandard system of phase field equations,
Contin. Mech. Thermodyn. {\bf 24} (2012), \pier{437-459}. 
}
		
\pier{\bibitem{CGS} P. Colli, G. Gilardi, J. Sprekels, 
Analysis and optimal boundary control of a nonstandard system of phase field equations, 
Milan J. Math. {\bf 80} (2012), 119-149.}

\pier{
\bibitem{CMR}
P. Colli, G. Marinoschi, E. Rocca, 
Sharp interface control in a Penrose-Fife model,
preprint arXiv:1403.4446 [math.AP] (2014), pp.~1-33.}

\betti{
\bibitem{DKS}
A. Damlamian, N. Kenmochi, N. Sato, 
Subdifferential operator approach to a class of nonlinear systems for 
Stefan problems with phase relaxation, Nonlinear Anal. {\bf 23} (1994), 115-142.
}

\betti{
\bibitem{EllZheng}
\pier{C.M. Elliott, S. Zheng,} 
Global existence and stability of solutions to the phase-field equations, 
in ``Free boundary problems'', Internat. Ser. Numer. Math., {\bf 95}, 46-58, Birkh\"auser
Verlag, Basel, (1990).}

\betti{
\bibitem{GraPetSch}
M. Grasselli, H. Petzeltov\'a, G. Schimperna,
Long time behavior of solutions to the Caginalp system with singular potential, 
Z. Anal. Anwend. {\bf 25} (2006), \pier{51-72}.
}

\betti{
\bibitem{HoffJiang}
K.-H.Hoffmann, L.S. Jiang, 
Optimal control of a phase field model for solidification,
Numer. Funct. Anal. Optim. \pier{{\bf 13} (1992), 11-27}. 
}

\betti{
\bibitem{HKKY}
K.-H. Hoffmann, N. Kenmochi, M. Kubo, N. Yamazaki, 
Optimal control problems for models of phase-field type with hysteresis of play operator,
Adv. Math. Sci. Appl. {\bf 17} (2007), \pier{305-336}.
}

\betti{
\bibitem{KenmNiez}
N. Kenmochi, M. Niezg\'odka,
Evolution systems of nonlinear variational inequalities arising phase change problems, 
\pier{Nonlinear} Anal. {\bf 22} (1994), 1163-1180.
}

\bibitem{LSU}
O.A. Lady\v zenskaja, V.A. Solonnikov, N.N. Ural'ceva:
``Linear and quasilinear equations of parabolic type'',
Trans. Amer. Math. Soc., {\bf 23},
Amer. Math. Soc., Providence, RI,
1968.

\betti{
\bibitem{Lau}
Ph. Lauren\c cot,
Long-time behaviour for a model of phase-field type, 
Proc. Roy. Soc. Edinburgh \pier{Sect.~A {\bf 126}} (1996), 167-185.
}
		
\pier{\bibitem{LK} C. Lefter, J. Sprekels, 
Optimal boundary control of a phase field 
system modeling nonisothermal phase transitions,
Adv. Math. Sci. Appl. {\bf 17} (2007), 181-194.}
		
\bibitem{Lions}
J.-L. Lions,
``\'Equations diff\'erentielles op\'erationnelles et probl\`emes aux limites'',
Grundlehren, Band~111,
Springer-Verlag, Berlin, 1961.

\betti{
\bibitem{Sch}
G. Schimperna, Abstract approach to evolution equations of phase field type and applications,
J. Differential Equations {\bf 164} (2000), 395-430.
}

\pier{\bibitem{SY} 
K. Shirakawa, N. Yamazaki, Optimal control problems of phase field 
system with total variation functional as the interfacial energy,
Adv. Differential Equations {\bf 18} (2013), 309-350.}

\bibitem{Simon}
J. Simon,
{Compact sets in the space $L^p(0,T; B)$},
\pier{Ann. Mat. Pura Appl.~(4)\/} 
{\bf 146} (1987), 65-96.

\betti{
\bibitem{SprZheng}
J. Sprekels, S. Zheng, 
Optimal control problems for a thermodynamically consistent model of phase-field type for phase transitions, Adv. Math. Sci. Appl. \pier{{\bf 1} (1992), 113-125}.
}

\End{thebibliography}

\End{document}

\bye